\documentclass[twocolumn]{article}
\usepackage[letterpaper,margin=1in]{geometry}
\usepackage{mathtools}
\usepackage[utf8]{inputenc}
\usepackage{xcolor}
\usepackage{cite}
\usepackage{hyperref}
\usepackage{mathscinet}
\usepackage{caption}
\usepackage{color,graphicx,enumerate,wrapfig,amssymb}
\usepackage{microtype}
\usepackage{stix2}

\newcommand{\R}{\mathbb{R}}

\makeatletter
\def\blfootnote{\gdef\@thefnmark{}\@footnotetext}
\makeatother

\makeatletter
\newcommand\HUGE{\@setfontsize\Huge{30}{40}}
\makeatother

\title{{\HUGE\bf Nina Nikolaevna
    Uraltseva}\\[2ex]\includegraphics[width=\textwidth]{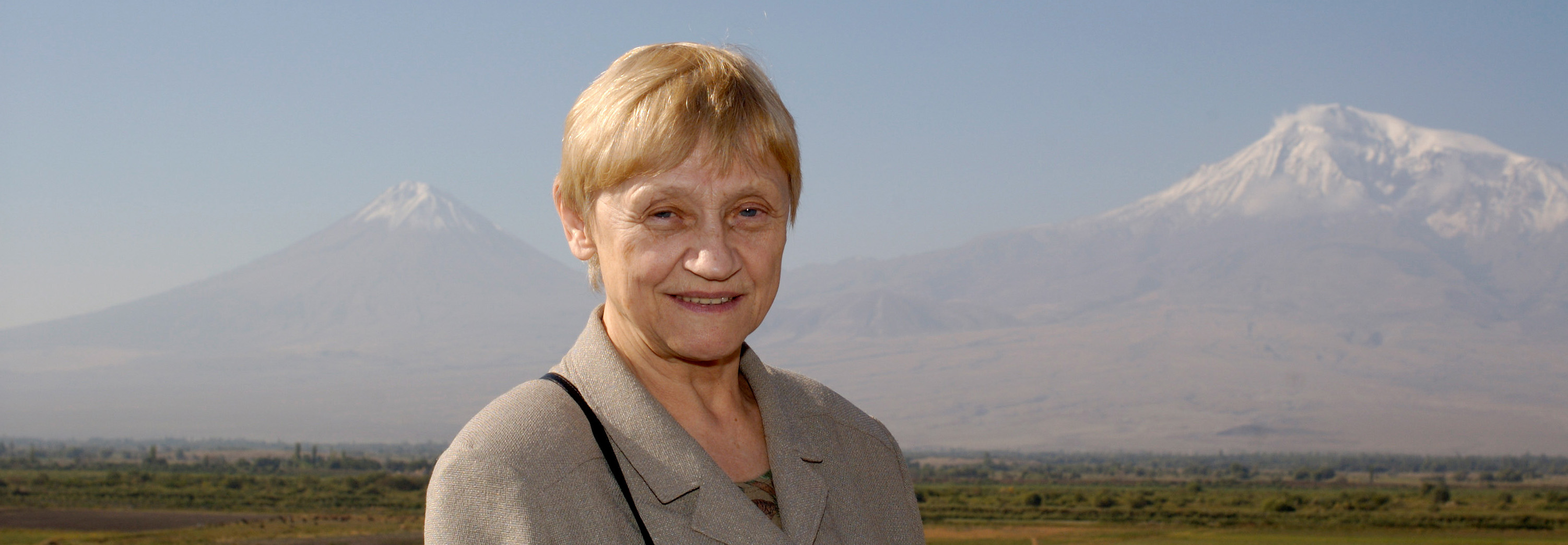}
}

\author{%
{\em\Large Darya Apushkinskaya}\thanks{S.M.\ Nikol'skii Mathematical Institute, Peoples' Friendship University of Russia (RUDN University), 117198 Moscow, Russia and Department of Mathematics and Computer Science, St. Petersburg State University, 199178 St. Petersburg, Russia;
 email: \href{mailto:apushkinskaya@gmail.com}{\tt apushkinskaya@gmail.com}} 
\and 
{\em\Large Arshak Petrosyan}\thanks{Department of Mathematics, Purdue University, West Lafayette, IN 47907, USA; email: \href{mailto:arshak@purdue.edu}{\tt arshak@purdue.edu}}
\and
{\em\Large Henrik Shahgholian}\thanks{Department of Mathematics, KTH Royal Institute of Technology, 10044 Stockholm, Sweden; email: \href{mailto:henriksh@kth.se}{\tt henriksh@kth.se}}%
}

\date{}

\begin{document}
\maketitle

\blfootnote{The opening photo is Nina Uraltseva with Mount Ararat in the background, Khor Virap, Armenia, 2004}

\noindent
Nina Nikolaevna Uraltseva was born on May 24, 1934 in Leningrad, USSR (currently St.\ Petersburg, Russia), to parents Nikolai Fedorovich Uraltsev, (an engineer) and  Lidiya Ivanovna Zmanovskaya (a school physics teacher). Nina Uraltseva was attracted to both mathematics and physics from the early stages of her life\footnote{Uraltseva's prematurely deceased younger brother (Igor Uraltsev) was a famous physicist, a specialist in epsilon spectroscopy in semiconductors. The Spin Optics Laboratory at St.\ Petersburg State University is named after him.}. She was a student at now famous school no.~239, then a school for girls, which later became specialized in mathematics and physics and produced many notable alumni.
Together with her friends, Nina Uraltseva initiated a mathematical study group at her school, under the supervision of  Mikhail Birman, then a student at the Faculty of Mathematics and Mechanics of Leningrad State University (LSU). In the higher grades of the school, she was actively involved in the Mathematical Circle at the Palace of Young Pioneers, guided by Ilya Bakelman, and became a two-time winner of the citywide mathematical olympiad.

Nina Uraltseva graduated from school in 1951 (with the highest distinction---a gold medal) and started her study at the Faculty of Physics of LSU. She was an active participant in an (undergraduate) student work-group founded by Olga Aleksandrovna Ladyzhenskaya, that gave her the opportunity to further deepening into the analysis of partial differential equations (PDEs). In 1956, she  graduated from the university and the same year married Gennady Lvovich Bir (a fellow student at the Faculty of Physics). The young couple  were soon blessed 
with a son (and the only child) Kolya.\footnote{Tragically, Kolya (Nikolai Uraltsev) passed away from a heart attack in 2013 (in Siegen, Germany). He was a renowned nuclear physicist, author of 120 papers published in the world's top scientific journals, most of them very well known internationally (with approximately 6000 references), and two of them are in the category of renowned. Kolya's son, Gennady Uraltsev, is currently a postdoctoral fellow at the University of Virginia, working in harmonic analysis.}

\begin{figure}[t]
\centering
\includegraphics[height=\columnwidth]{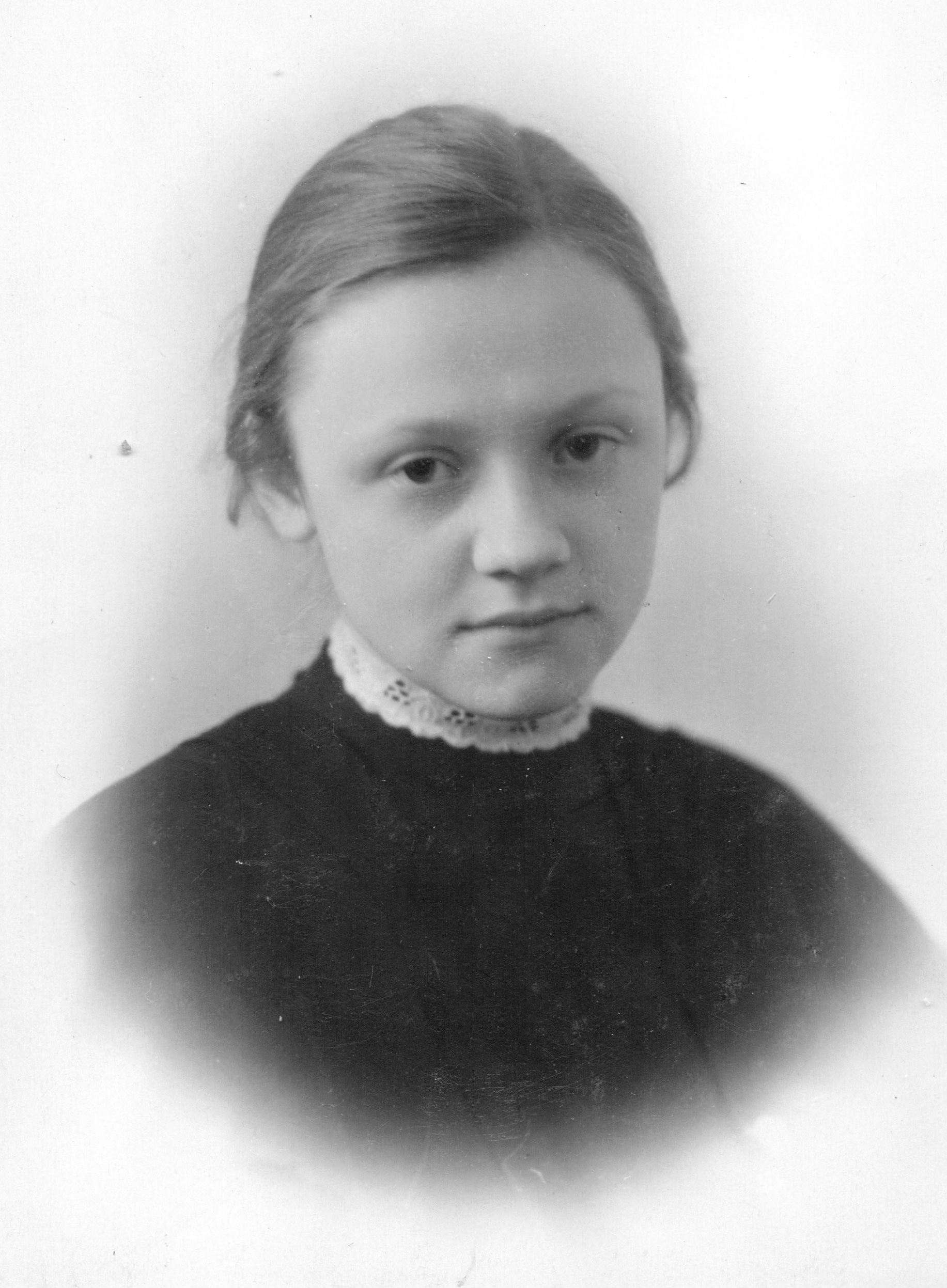}
\caption{Nina Uraltseva in a  schoolgirl uniform, Leningrad, 1951.}
\label{fig:1951}
\end{figure}

During  her graduate years,  Uraltseva continued to be supervised by Olga Ladyzhenskaya. This mentorship transformed to a lifelong productive collaboration and warm friendship  until 2004, when Olga Ladyzhenskaya passed away.

Nina Uraltseva defended her Candidate of Science\footnote{Equivalent of PhD in many countries.} thesis entitled  ``Regularity of solutions to multidimensional quasilinear equations and variational problems'' in 1960. Four years later, she  became a Doctor of Science\footnote{Equivalent of Habilitation in many European countries.} with a thesis ``Boundary-value problems for quasilinear elliptic equations and systems of second order.'' 
Since 1959, she has been a member of the Chair of Mathematical Physics at the Faculty of Mathematics and Mechanics of LSU (currently  St.\ Petersburg State University), where she became a Full Professor in 1968  and served as the head of the chair since 1974.

For her fundamental contributions to the theory of partial differential equations in the 1960s, Nina Uraltseva (jointly with Olga Ladyzhenskaya) was awarded the Chebyshev Prize of the Academy of Sciences of the USSR (1966) as well as one of the highest honors of the USSR, the USSR State Prize (1969). 

Throughout her career, Nina Uraltseva has been an invited speaker in many meetings and conferences, including the International Congress of Mathematicians in 1970 and 1986. In 2005, she was chosen as the Lecturer of the European Mathematical Society.

Nina Uraltseva's mathematical achievements are highly regarded throughout the world, and have been acknowledged by various  awards, such as  the titles of Honorary Scientist of the Russian Federation in 2000, Honorary Professor of St.\ Petersburg State University in 2003, and Honorary Doctor of KTH Royal Institute of Technology, Stockholm, Sweden, in 2006. In the same year, in recognition of her academic record, she received the Alexander von Humboldt Research Award. In 2017, the Government of St.\ Petersburg recognized her recent research by its Chebyshev Award.

Nina Uraltseva's interests are not limited to scientific activities only. In her youth, she used to be  a very good basketball player and   an active member of the university basketball team. She enjoyed hiking in the mountains, canoeing, and car driving. In the 1980s, Nina took part in five archaeological expeditions in the north of Russia (the Kola Peninsula and the Kotlas area) and excavated Paleolithic ceramics. She is also a passionate lover of classical music and a regular visitor at philharmonic concerts.

\section*{Mathematical Contributions}

Nina Uraltseva has made lasting contributions to mathematics with her pioneering work in various directions in analysis and PDEs and the development of elegant and sophisticated analytical techniques. She is most renowned for her early work on linear and quasilinear equations of elliptic and parabolic type in collaboration with Olga Ladyzhenskaya, which is the category of classics, but her contributions to the other areas such as degenerate and geometric equations, variational inequalities, and free boundaries are equally deep and significant. Below, we 
summarize  Nina Uraltseva's work with some details on selected results.

\subsection{Linear and Quasilinear Equations} 

\subsubsection{Hilbert's 19th and 20th problems}

The first three decades of 
Nina Uraltseva's mathematical career were devoted to the theory of linear and quasilinear PDEs of elliptic and parabolic type. Her first round of works in the 1960s, mostly in collaboration with Olga Ladyzhenskaya, was related to Hilbert's 19th and 20th problem on the existence and regularity of the minimizers of the energy integrals
$$
I(u)=\int_{\Omega} F(x,u,\nabla u)dx,
$$
where $F(x,u,p)$ is a smooth function of its arguments and $\Omega$ is a bounded domain in $\R^n$, $n\geq 2$. In her Candidate of Science thesis, based on work \cite{MR0126742},\footnote{In those years, it was quite unusual to base the Candidate of Science thesis on just a single paper and some of the committee members voiced their concerns. 
However, Olga Ladyzhenskaya objected decisively that it depends on the quality of the paper.} 
Nina Uraltseva has shown that under the assumption that $F$ is $C^{2,\alpha}$ and satisfies the uniform ellipticity condition
$$
F_{p_ip_j}\xi_i\xi_j\geq m |\xi|^2,\quad m>0,
$$
the minimizers $u$ are $C^{2,\alpha}$ locally in $\Omega$ (i.e., on compact subdomains of $\Omega$), provided they are  Lipschitz. (It has to be mentioned here that the Lipschitz regularity of the minimizers was known from the earlier works of Ladyzhenskaya under natural growth condition on $F$ and its partial derivatives.) Uraltseva has also shown that $C^{2,\alpha}$ regularity extends up to the boundary $\partial\Omega$ under the natural requirement that both $\partial\Omega$ and $u|_{\partial\Omega}$ are $C^{2,\alpha}$. This generalized the results of  Morrey in dimension $n=2$ to higher dimensions.  

\begin{figure}[t!]
\centering
  \includegraphics[width=\columnwidth]{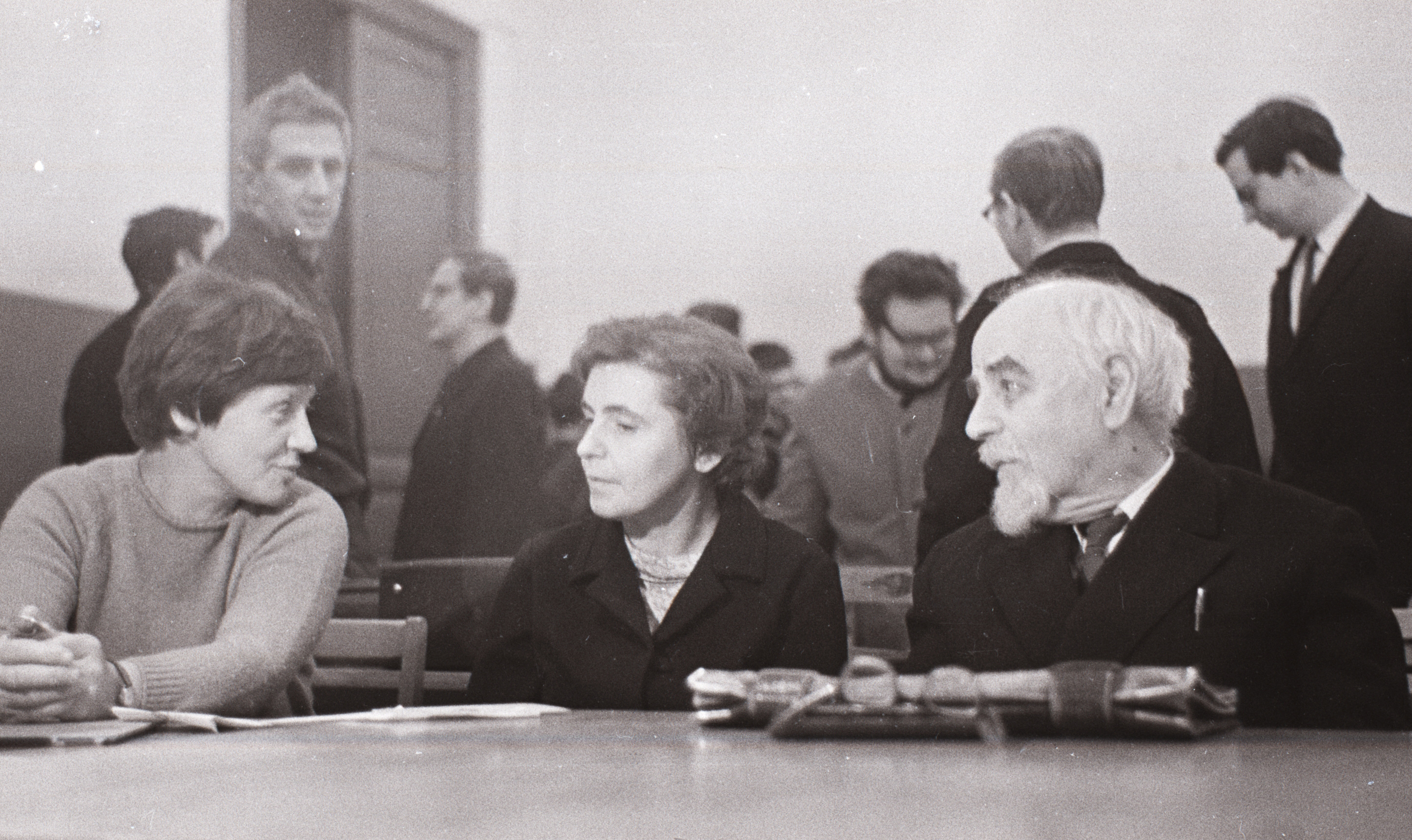}
\caption{Left to right: Nina Uraltseva, Olga Ladyzhenskaya, and Vladimir Smirnov in a seminar on mathematical physics, Leningrad, 1968.}
\label{fig:seminar}
\end{figure}

Uraltseva's proof was based on a deep extension of the ideas of De Giorgi for the solutions of uniformly elliptic equations in divergence form with bounded measurable coefficients, which were applicable only to the integrands of the form $F(\nabla u)$. 
In particular, one of the essential steps was to establish that $v=\pm u_{x_i}$, $i=1,\dots,n$, which are assumed to be bounded,
satisfy the energy inequalities
\begin{equation}\label{eq:energy-ineq}
\int_{A_{k,\rho}}|\nabla v|^2\zeta^2\leq C\int_{A_{k,\rho}}(v-k)^2|\nabla \zeta|^2+C|A_{k,\rho}|
\end{equation}
for all $|k|\leq M$, where $A_{k,\rho}$ is intersections of $\{v>k\}$ with the ball $B_\rho(x^0)\Subset \Omega$, $\zeta$ is a cutoff function, and $M$ is a bound for $\max|\nabla u|$.

Using similar ideas, Uraltseva was able to deduce the existence and regularity of solutions for the class of quasilinear uniformly elliptic equations in divergence form,
\begin{equation}\label{eq:div}
\partial_{x_i}(a_{i}(x,u,\nabla u))+a(x,u,\nabla u)=0,
\end{equation}
under natural growth conditions on $a_i(x,u,p)$, $a$ and some of their partial derivatives, which were mainly needed for proving the bounds on $\max |\nabla u|$. These results were further refined in the joint works with Olga Ladyzhenskaya \cite{MR0141891,MR0141874,MR0149075,MR149076,MR0150447} as well as in \cite{MR0142886}, for the case of Neumann-type boundary conditions. The latter paper also contained similar results for certain quasilinear diagonal systems (important, e.g., for the applications in harmonic maps).

Quasilinear uniformly elliptic equations in nondivergence form,
\begin{equation}\label{eq:nondiv}
a_{ij}(x,u,\nabla u)\,u_{x_ix_j}+a(x,u,\nabla u)=0,
\end{equation}
were trickier to treat, but already in \cite{MR0126742},
Uraltseva found a key: quadratic growth of $a(x,u,p)$ in the $p$-variable,
$$
|a(x,u,p)|\leq \mu(1+|p|)^2,
$$
along with the corresponding conditions on the partial derivatives of $a$ and $a_{ij}$ in their variables. In \cite{MR0140817}, Uraltseva proved $C^{1,\alpha}$ a priori bounds for solutions of \eqref{eq:nondiv}, as wells as for diagonal systems of similar type.

The results in the elliptic case were further extended to the parabolic case (including systems) in a series of works of Ladyzheskaya and Uraltseva \cite{MR0147786,MR0181837}. 

\begin{figure}
    \centering
    \begin{picture}(201,150)
    \put(0,0){%
    \includegraphics[height=150pt]{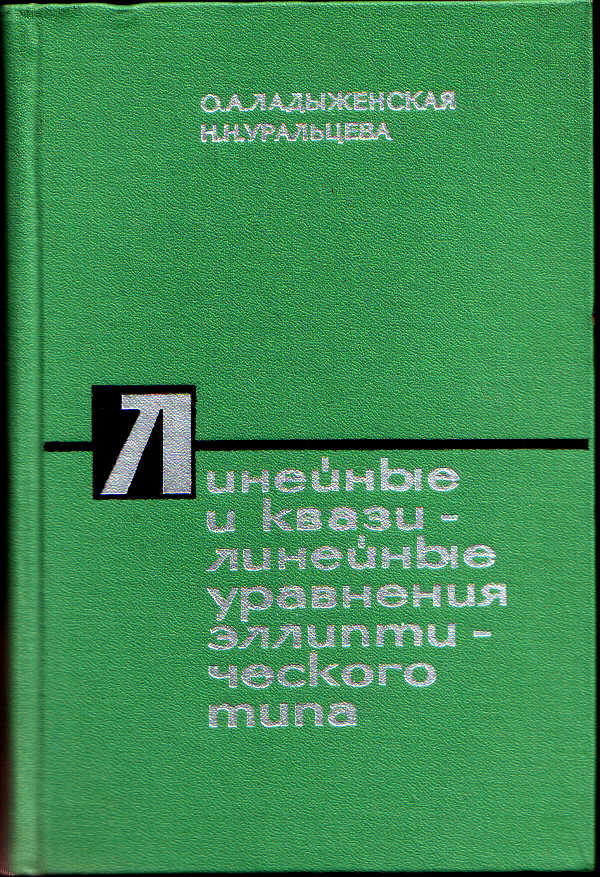}%
    \includegraphics[height=150pt]{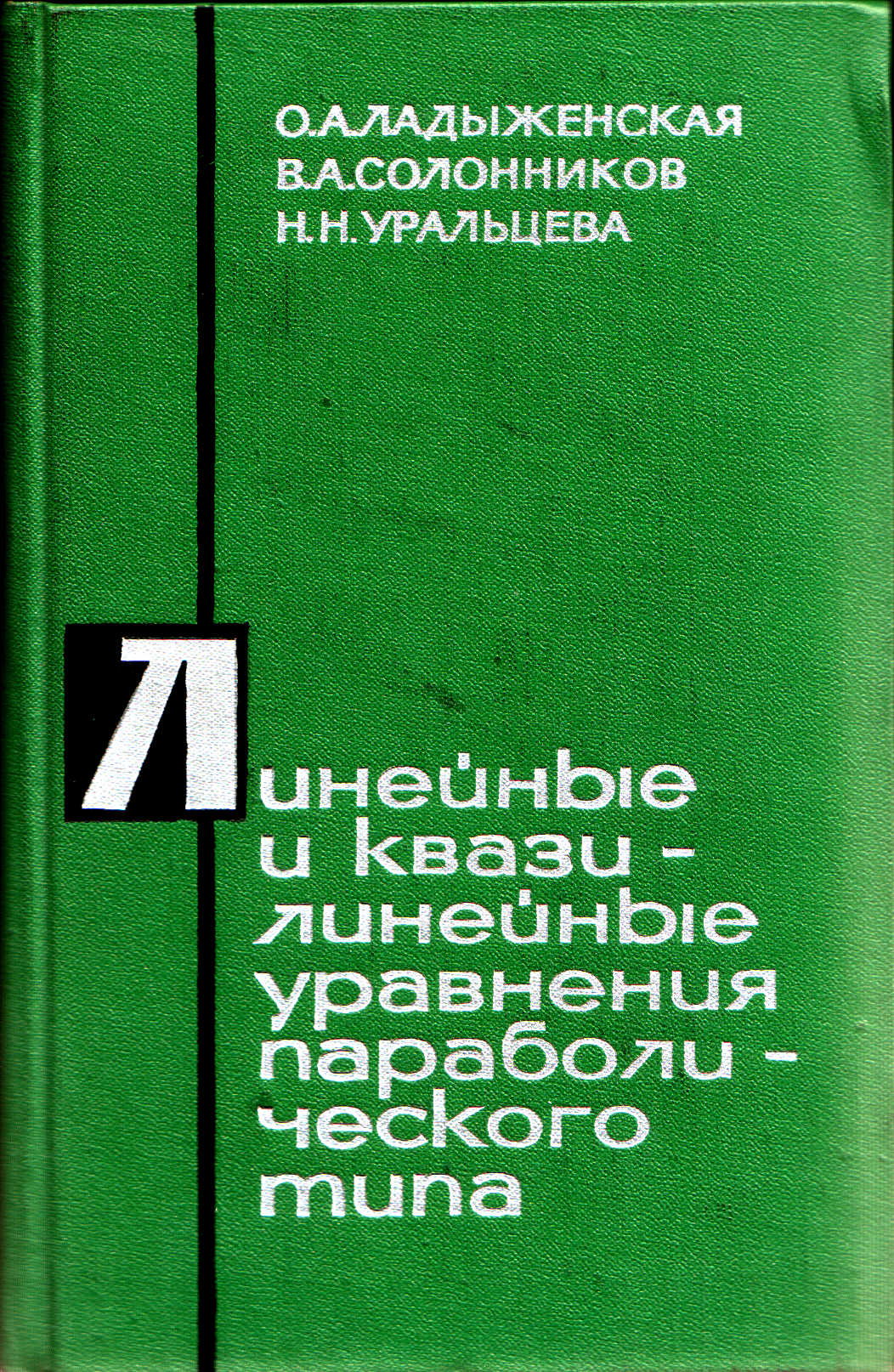}}
    \end{picture}
    \caption{The famous books: the iconic green Russian editions of the Elliptic (2nd ed., 1973) and Parabolic (1967) versions of Uraltseva's books with Ladyzhenskaya and Solonnikov.}
    \label{fig:books}
\end{figure}

This extensive research, that went far beyond the original scope of Hilbert's 19th and 20th problems, was summarized in two monographs, \emph{Linear and Quasilinear Equations of Elliptic Type} \cite{MR0211073} (substantially enhanced in the 2nd edition \cite{MR0509265}) and \emph{Linear and Quasilinear Equations of Parabolic Type} \cite{MR0241821}, written in collaboration with Vsevolod Solonnikov; see Figure~\ref{fig:books}. The monographs became instant classics and were translated to English \cite{MR0244627,MR0241822} and other languages and have been extensively used for generations of mathematicians working in elliptic and parabolic PDEs and remain so to this date.

\subsubsection{Equations with unbounded coefficients}

In a series of papers in 1979--1985, summarized in her talk at the International Congress of Mathematicians in Berkeley, CA, 1986 \cite{MR934318} and a survey paper with Ladyzhenskaya \cite{MR878325}, Uraltseva and collaborators have studied uniformly elliptic quasilinear equations of nondivergence type \eqref{eq:nondiv} and their parabolic counterparts, when $a$ and the first derivatives of $a_{ij}$ are possibly unbounded. The typical conditions read
$$
|a(x,u,p)|\leq \mu |p|^2+b(x)|p|+\Phi(x),
$$
where $\mu$ is a constant and $b,\Phi\in L^q(\Omega)$, $q>n$. Uraltseva and collaborators were able to establish the existence and up to the boundary $C^{1,\alpha}$ regularity of $W^{2,n}$ strong solutions of the problem, vanishing on $\partial\Omega$ (provided the latter is sufficiently regular). The proofs were based on the extension of methods of Ladyzhenskaya and Uraltseva already in their books \cite{MR0244627,MR0241822}, as well as those of Krylov and Safonov using the Aleksandrov-Bakelman-Pucci (ABP) estimate, in the elliptic case, and a parabolic version of the ABP estimate due to Nazarov and Uraltseva \cite{MR821477}, in the parabolic case.

Most recent results of Nina Uraltseva in this direction are in the joint work with Alexander Nazarov \cite{MR2760150} on the linear equations in divergence form,
$$
\partial_{x_i}(a_{ij}(x)u_{x_j})+b_i(x)u_{x_i}=0\quad\text{in }\Omega,
$$
and their parabolic counterparts. Their goal was to find conditions on the lower order coefficients $\mathbf{b}=(b_1,\dots,b_n)$ that guarantee the validity of classical results such as the strong maximum principle, Harnack's inequality, and Liouville's theorem. It was shown by Trudinger \cite{MR369884} that such results hold when $\mathbf{b}\in L^q$, $q>n$. Motivated by applications in fluid dynamics, in one of their theorems, Nazarov and Uraltseva showed that under the additional assumption,
\begin{equation}
\label{eq:div-free}
\operatorname{div}\mathbf{b}=0,
\end{equation}
the condition on $\mathbf{b}$ can be relaxed to being in the Morrey space
$$
\sup_{B_r(x^0)\subset\Omega}r^{q-n}\int_{B_r(x^0)}|\mathbf{b}|^q<\infty
$$
for some $n/2<q\leq n$. In the borderline case $q=n$, the Morrey space above is locally the same as $L^n$. Remarkably, in that case the divergence free condition \eqref{eq:div-free} on $\mathbf{b}$ can be dropped when $n\geq 3$, i.e., $\mathbf{b}\in L^n_{\mathrm{loc}}$ alone is sufficient to have the classical theorems; moreover, this result is optimal. In dimension $n=2$, to drop \eqref{eq:div-free} one needs a stronger condition $\mathbf{b}\ln^{1/2}(1+|\mathbf{b}|)\in L^2_{\mathrm{loc}}$.

\subsection{Nonuniformly Elliptic and Parabolic Equations}

\subsubsection{Degenerate Equations} Nina Ural\-tse\-va has also made a pioneering  work on the regularity theory for degenerate quasilinear equations. A particular result in this direction is her 1968 proof \cite{MR0244628} of the $C^{1,\alpha}$ regularity of $p$-harmonic functions, $p>2$, which are the weak solutions of the $p$-Laplace equation
 \begin{equation}\label{p-lap}
 \operatorname{div}(|\nabla u|^{p-2}\nabla u)=0\quad\text{in }\Omega,
 \end{equation}
or, equivalently, are the minimizers of the energy functional
 $$
 \int_{\Omega}|\nabla u|^pdx.
 $$
The difficulty here lies in the fact that the $p$-Laplace equation \eqref{p-lap} degenerates at the points where the gradient vanishes and that the solutions are not generally twice differentiable in the Sobolev sense. As stated in her paper, this problem was posed to Nina Uraltseva by Yurii Reshetnyak in relation with the study of quasiconformal mappings in higher dimensions.
 
Uraltseva has obtained the $C^{1,\alpha}$ regularity of $p$-harmonic functions as an application of the H\"older regularity of the solutions of the degenerate quasilinear diagonal systems 
$$
\partial_{x_i}(a_{ij}(x,\mathbf{u})\,\mathbf{u}_{x_j})=\mathbf{0},
$$
with scalar coefficients $a_{ij}$ satisfying the degenerate ellipticity condition
$$
\nu(|\mathbf{u}|)|\xi|^2\leq a_{ij}(x,\mathbf{u})\xi_i\xi_j\leq \mu \nu(|\mathbf{u}|)|\xi|^2,
$$
with $\mu\geq 1$ and a nonnegative increasing function $\nu(\tau)$ satisfying $\nu(\lambda\tau)\leq \lambda^s\nu(\tau)$ for $\lambda\geq 1$ and $s>0$.

Unfortunately, despite the utmost importance of this result, Nina Uraltseva's proof remained unknown outside of the Soviet Union. In 1977, nine years later, it was independently reproved by Karen Uhlenbeck \cite{MR474389}. Other proofs were given by Craig Evans \cite{MR672713}, John Lewis \cite{MR721568}, who extended the range of exponents to $1<p\leq 2$, Di Benedetto \cite{MR709038} and Tolksdorf \cite{MR727034}, who both extended it to the case of general degenerate quasilinear equations in divergence form.  

Another work in this area that has gained the status of classic is the paper of Nina Uraltseva and Anarkul Urdaletova \cite{MR725829}, where they proved uniform gradient estimates for bounded solutions of anisotropic degenerate equations,
$$
\partial_{x_i}(a_i(x,u_{x_i}))+a(x,u,\nabla u)=0\quad\text{in }\Omega,
$$
under ellipticity, growth, and monotonicity conditions on the coefficients. Their results were applicable to the minimizers of the energy functional
$$
\int_{\Omega} \sum_{i=1}^n|u_{x_i}|^{m_i}+f(x,u),
$$
with the exponents $m_1,\dots,m_n$ satisfying $m_i>3$, $2m_i>m_0$, $i=1,\dots,n$, $m_0=\max\{m_i\}$, under the monotonicity condition $f_u(x,u)\geq 0$. This was the very first paper to prove regularity results for degenerate quasilinear equations with nonstandard growth, which appeared first in the 1980s, motivated by applications in elasticity and material science, and continue to be the subject of extensive research today \cite{MR4258810}. Major contributions in this direction have been made by Paolo Marcellini \cite{MR868523,MR969900,MR1094446} and many others.

\subsubsection{Geometric Equations}
\label{sec:geom-eq}

In \cite{MR265745}, Ladyzhenskaya and Uraltseva developed a method of local a priori estimates for nonuniformly elliptic and parabolic equations, including the equations of minimal surface type,
$$
\operatorname{div}\frac{\nabla u}{\sqrt{1+|\nabla u|^2}}=a(x,u,\nabla u)\quad\text{in }\Omega.
$$
A particular case with $a(x,u,\nabla u)=\kappa u$, $\kappa>0$, together with the Neumann-type condition 
${\partial_\nu u}/{\sqrt{1+|\nabla u|^2}}=\varkappa$
on $\partial\Omega$,  $|\varkappa|<1$, 
is known as the capillarity problem. The boundary estimates, as well as the existence of classical solutions for such problems were obtained in \cite{MR0638359,MR0638360,MR660089,MR0364860}. Remarkably, the estimates in the last two papers required only the smoothness of the domain $\Omega$, but not its convexity.

In the 1990s, in a series of joint works \cite{MR1193345,MR1246345,MR1334142,MR1483640} with Vladimir Oliker, Nina Uraltseva studied the evolution of surfaces $S(t)$ given as graphs $u = u(x, t)$ over a bounded domain $\Omega\subset \R^n$ with the speed depending on the mean curvature of $S(t)$ under the condition that the boundary of the surface $S(t)$ is fixed. More precisely, they considered a parabolic PDE of the type
$$
u_t=\sqrt{1+|\nabla u|^2}\,\operatorname{div}\frac{\nabla u}{\sqrt{1+|\nabla u|^2}}\quad\text{in }\Omega\times(0,\infty)
$$
with the boundary condition $u(x,t)=\phi(x)$ on $\partial\Omega\times(0,\infty)$ and initial condition $u(x,0)=u_0(x)$. Even in the stationary case, when this problem is the Dirichlet problem for the mean curvature equation, the existence of up to the boundary classical solutions requires a geometric condition on the domain $\Omega$, namely, the nonnegativity of the mean curvature of $\partial\Omega$. For such domains, Huisken \cite{MR983300} has shown the existence of the classical solutions of the evolution problem and proved that the surfaces $S(t)$ converge to a classical minimal surface $S$ as $t\to\infty$. Oliker and Uraltseva have studied this problem with no geometric conditions on the domain $\Omega$. For this purpose, they introduced a notion of a generalized solution to the parabolic problem (as a limit of regularized problems). They have proved its existence and convergence $u(\cdot,t)\to \Phi$ as $t\to\infty$ to a generalized solution $\Phi$ of the stationary problem, in the sense that $\Phi$ minimizes the area functional
$$
\int_\Omega \sqrt{1+|\nabla u|^2}+\int_{\partial\Omega} |u-\phi|
$$
among all competitors in $W^{1,1}(\Omega)$. Such minimizer $\Phi$ is unique, but may differ from the Dirichlet data $\phi$ on the ``bad'' part of the boundary where the mean curvature is negative. The study of the behavior of the minimizer near the ``contact points'' on the boundary where $\Phi|_{\partial \Omega}$ ``detaches'' from $\phi$ later served as one of Uraltseva's motivations for studying the touch between free and fixed boundaries; see Section~\ref{sec:touch}.

\subsection{Variational Inequalities}

Another area in which Nina Uraltseva has made significant contributions is variational inequalities, including variational problems with convex constraints that often exhibit a priori unknown sets known as free boundaries. An important example is the Signorini problem from elasticity, which describes equilibrium configurations of an elastic body resting on a rigid frictionless surface.

In a series of papers \cite{MR775926,MR1374941,MR841488,MR928154,MR964910,MR0313623,MR508511,MR869237,MR860572}, partially with  Arina Arkhipova,
Nina Uraltseva studied elliptic and parabolic variational inequalities
with unilateral and bilateral boundary constraints, known as the boundary obstacle problems, which can be viewed as scalar versions of the Signorini problem. Ultimately, these results played a fundamental role in Schumann's proof \cite{MR986184} of the
 $C^{1,\alpha}$ regularity for the solution of the Signorini problem in the vectorial case.

Below, we give a more detailed description of some of her most impactful results in this direction.

\subsubsection{Problems with unilateral constraints}

\begin{figure}[t]
\centering
\begin{picture}(123,150)(0,0)
\put(0,0){\includegraphics[height=150pt]{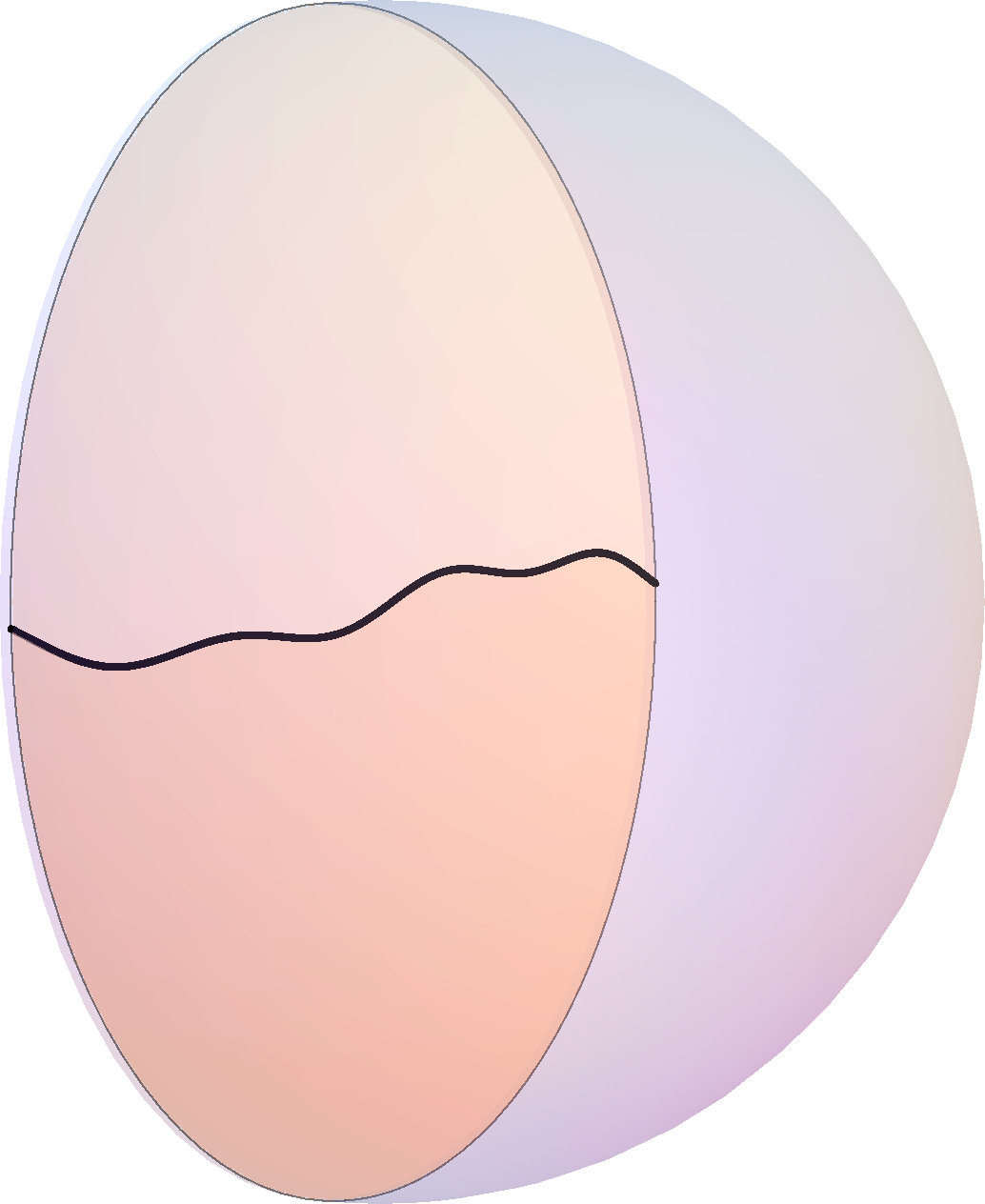}}			
\put(95,75){\footnotesize $\Omega$}	
\put(22,128){\footnotesize $S$}
\put(26,48){\footnotesize $u>\psi$}
\put(24,33){\footnotesize $\partial_\nu^A u=0$}
\put(27,106){\footnotesize $u=\psi$}
\put(25,91){\footnotesize $\partial_\nu^A u\geq 0$}
\put(68,68){\footnotesize $\Gamma$}
\end{picture}
\caption{Boundary obstacle problem.}\label{fig:bdry-obst}
\end{figure}
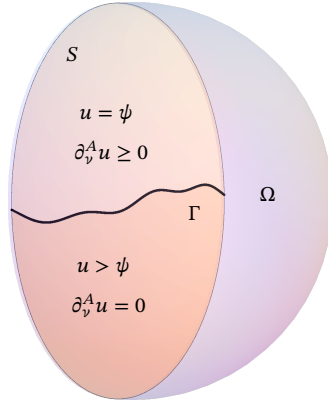

Let
$\Omega\subset\R^n$, $n\geq 2$, be a bounded domain with a smooth
boundary and $S$ a relatively open nonempty subset of $\partial
\Omega$. Suppose we are also given two functions $\psi,g\in
W^{1,2}(\Omega)$ satisfying $g\geq \psi$ on $S$ (in the
sense of traces). Consider then a closed convex subset
$\mathfrak{K}\subset W^{1,2}(\Omega)$ defined by
$$
\mathfrak{K}\coloneq\{v\in W^{1,2}(\Omega): \text{$v\geq \psi$ on
  $S$, $v=g$ on $\partial\Omega\setminus S$}\}.
$$
In other words, $\mathfrak{K}$ consists of functions that need to stay above $\psi$, called
a boundary (or thin) obstacle, on
$S$ and equal to $g$ on
$\partial\Omega\setminus S$. Then, one wants to find
$u\in\mathfrak{K}$ that minimizes the generalized Dirichlet energy
$$
J(v)=\int_\Omega a_{ij}(x)v_{x_j}v_{x_i}+2f(x) u,
$$
where $a_{ij}(x)$ are uniformly elliptic coefficients and $f$ is a
certain function. Equivalently, the minimizer $u$ satisfies the
variational inequality
\begin{multline*}
u\in\mathfrak{K},\quad\int_\Omega  a_{ij}(x)u_{x_j}(v-u)_{x_i}\\+f(x)(v-u)\geq
0\quad\text{for any }v\in\mathfrak{K}.
\end{multline*}
In turn, it is equivalent to the following boundary value
problem
\begin{align*}
  \partial_{x_i}(a_{ij}(x)u_{x_j})=f(x)&\quad\text{in }\Omega,\\
  u=g&\quad\text{on }\partial\Omega\setminus S,\\
  u\geq \psi,\ \partial_\nu^A u\geq 0,\ (u-\psi)\partial_\nu^A u=0&\quad\text{on }S,
\end{align*}
to be understood in the appropriate weak sense,
where $\partial_\nu^A u\coloneq a_{ij}(x)\nu_j u_{x_j}$ is the conormal derivative of $u$ on $\partial \Omega$, with
$\nu=(\nu_1,\dots,\nu_n)$ being the outward unit normal. 
The conditions on
$S$ are known as the 
Signorini complementarity conditions
and are remarkable because they imply that
$$
\text{either $u=\psi$ or $\partial_\nu^A u=0$ on $S$},
$$
yet the exact sets where the first or the second equality holds are
unknown. The interface $\Gamma$ between these sets in $S$ is called free boundary  (see Figure~\ref{fig:bdry-obst}). The study of the free boundary is one of the main objectives in such problems (see  Section~\ref{sec:free-bdry-prob} for Uraltseva's contributions in that direction), yet the regularity of the solutions $u$ is a challenging problem by itself and is often an important step towards the study of the free boundary.

One of the theorems of Nina Uraltseva \cite{MR860572,MR775926} states that when
\begin{align*}
a_{ij}&\in W^{1,q}(\Omega),\\
\psi&\in W^{1,2}(\Omega)\cap
W^{2,q}_{\mathrm{loc}}(\Omega\cup S),\\
f&\in
L^q(\Omega),
\end{align*}
for some $q>n$, then
$$
u\in C^{1,\alpha}_{\mathrm{loc}}(\Omega\cup S),
$$ 
with a universal exponent $\alpha\in(0,1)$. Prior to this result, similar conclusion was known only under higher regularity assumptions on the coefficients and the obstacle in the works of Caffarelli \cite{MR542512} and Kinderlehrer \cite{MR620584}.
The lower regularity assumptions in Uraltseva's result,  particularly on the obstacle $\psi$, were instrumental in Schumann's proof of the corresponding result in the vectorial case \cite{MR986184}. 
The parabolic counterpart of Uraltseva's theorem, with similar assumptions on the coefficients and the obstacle was established later in a joint work of Arkhipova and Uraltseva \cite{MR1374941}.

The idea of Uraltseva's proof is based on an interplay between De Giorgi-type energy inequalities and the Signorini complementarity condition. Locally, near $x^0\in S$, one can assume that $S=\{x_n=0\}$ and $\psi=0$. First,
working with the regularized problem, one can establish that for any partial derivative $v=\pm u_{x_i}$, $i=1,\dots,n$, there holds an energy inequality (similar to \eqref{eq:energy-ineq} in the unconstrained case)
$$
\int_{A_{k,\rho}}|\nabla v|^2\zeta^2\leq C\int_{A_{k,\rho}}(v-k)^2|\nabla\zeta|^2+C_0|A_{k,\rho}|^{1-2/q}
$$
for any $k>0$, $0<\rho<\rho_0$, and a cutoff function $\zeta$ in $B_\rho(x^0)$, where $A_{k,\rho}=\{v>k\}\cap B_\rho(x^0)\cap\Omega$.
Next, one observes that as a consequence of the Signorini complementarity conditions, one has
$$
u_{x_i}u_{x_n}=0\quad \text{on }\{x_n=0\}\cap B_\rho(x_0)
$$
for all $i=1,\dots,n-1$ and hence either the normal derivative $v=u_{x_n}$ or all tangential derivatives $v=u_{x_i}$, $i=1,\dots,n-1$, vanish at least on half of $\{x_n=0\}\cap B_\rho(x^0)$ (by measure). This allows to apply Poincare's inequality in one of the steps and obtain a geometric improvement of the Dirichlet energy for $v$ going from radius $\rho$ to $\rho/2$. By iteration, this gives that either
\begin{align}
\label{eq:est-tang}\sum_{i=1}^{n-1}\int_{\Omega\cap B_\rho(x^0)}|\nabla u_{x_i}|^2&\leq C\rho^{n-2+2\alpha}\quad\text{or}\\
\label{eq:est-norm}
\int_{\Omega\cap B_\rho(x^0)}|\nabla u_{x_n}|^2&\leq C\rho^{n-2+2\alpha}
\end{align}
holds, with $C$ depending on the distance from $x^0$ to $\partial\Omega\setminus S$. However, using the PDE satisfied by $u$, it is easy to see that \eqref{eq:est-tang} implies \eqref{eq:est-norm}, and hence \eqref{eq:est-norm} always holds. From there, the $C^{1,\alpha}$ regularity of $u$ follows by standard results for the solutions of the Neumann problem.

\subsubsection{Diagonal systems}
The results described above were extended, in joint works with Arina Arkhipova \cite{MR841488,MR964910}, to the problem with two obstacles $\psi_-\leq \psi_+$ on $S$, that corresponds to the constraint set
\begin{multline*}
\mathfrak{K}=\{v\in W^{1,2}(\Omega): \text{$\psi_-\leq u \leq \psi_+$ on $S$},\\\text{$u=g$ on $\partial\Omega\setminus S$}\}.
\end{multline*}
While substantial difficulties arise near the set where
$\psi_-=\psi_+$, the results are as strong as in the case of a single obstacle. In their further work, Arkhipova and Uraltseva \cite{MR869237,MR928154} studied related problems for quasilinear elliptic systems with diagonal principal part.
To describe their results, let $V=W^{1,2}(\Omega;\R^N)\cap L^\infty(\Omega;\R^N)$ and
$$
\mathfrak{K}=\{\mathbf{u}\in V: \text{$\mathbf{u}(x)\in K(x)$ for every $x\in \partial\Omega$}\},
$$
where $K(x)$ are given convex subset of $\R^N$ for every $x\in \partial\Omega$.
Then consider the variational inequality of the type
\begin{multline*}
   \mathbf{u}\in\mathfrak{K},\quad \int_{\Omega} \left(a_{ij}(x,\mathbf{u}) \mathbf{u}_{x_j}+\mathbf{b}_{i}(x,\mathbf{u})\right)(\mathbf{v}-\mathbf{u})_{x_i}\\
    +\mathbf{f}(x,\mathbf{u},\nabla \mathbf{u})(\mathbf{v}-\mathbf{u})\geq \mathbf{0}\\\text{for any }\mathbf{v}\in\mathfrak{K},
\end{multline*}
where $a_{ij}$ are scalar uniformly elliptic coefficients, $\mathbf{b}_i$ and $\mathbf{f}$ are $N$-dimensional vector functions and $\mathbf{f}(x,\mathbf{u},\mathbf{p})$ grows at most quadratically in $\mathbf{p}$. We note that the problem with two obstacles $\psi_-\leq\psi_+$ on $\partial\Omega$ fits into this framework with $N=1$ and $K(x)=[\psi_-(x),\psi_+(x)]$. Assume now that the convex sets $K(x)$ are of the form
$$
K(x)=T(x)K_0+\mathbf{g}(x),
$$
where $K_0$ is a convex set in $\R^N$ with a nonempty interior and a
smooth ($C^2$) boundary, $T(x)$ is an orthogonal $N\times N$ matrix,
and $\mathbf{g}(x)$ is an $N$-dimensional vector. A theorem of
Arkhipova and Uraltseva \cite{MR928154} then states that when the
entries of $T$ and $\mathbf{g}$ are extended to $W^{2,q}$ functions in
$\Omega$, $q>n$, $a_{ij}(\cdot,\mathbf{u})$ and $\mathbf{b}_i(\cdot,\mathbf{u})$ are in $W^{1,q}(\Omega)$,
uniformly in $\mathbf{u}$ and have at most linear growth in $\mathbf{u}$, and
$\mathbf{f}$ has at most quadratic growth in $\mathbf{p}$, then
$$
\mathbf{u}\in C^{1,\alpha}_{\text{loc}}(\Omega\cup S),
$$
provided $\mathbf{u}$ is H\"older continuous in $\Omega\cup S$. The H\"older continuity assumption on $\mathbf{u}$ can be replaced by a bound on the oscillation in $\Omega$ and a local uniqueness of the solutions, which is also necessary for the continuity of the solutions of the nonlinear systems of the type  
$$
\partial_{x_i}(a_{ij}(x,\mathbf{u})\mathbf{u}_{x_j}+\mathbf{b}_i(x,\mathbf{u}))+\mathbf{f}(x,\mathbf{u},\nabla\mathbf{u})=\mathbf{0}\quad \text{in }\Omega.
$$
For a more complete overview of Uraltseva's results on variational inequalities, we refer to her own survey paper \cite{MR933999}.


\subsection{Free Boundary Problems}
\label{sec:free-bdry-prob}

In the last 25 years, Uraltseva’s work has dealt with regularity issues arising in free boundary problems. She has developed powerful techniques, which have led to proving the optimal regularity results for solutions and for free boundaries. She has systematically studied how the free boundaries approach the fixed boundaries \cite{MR1950478}, and has developed tools to study free boundary problems for weakly coupled systems \cite{MR3350233}, as well as two-phase problems \cite{MR2340105}. The graduate textbook \emph{Regularity of Free Boundaries in Obstacle-Type Problems} \cite{MR2962060}, written in collaboration with two of us, contains these and related results.

Some of Uraltseva's major contributions (results, approaches) in free boundary problems are addressed below in more detail.

\begin{figure}
  \centering
\begin{picture}(206,100)
\put(0,10){\includegraphics[height=90pt]{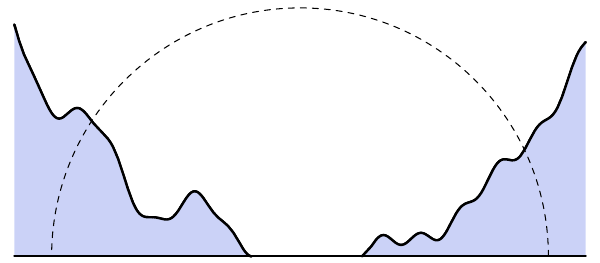}}
\put(10,0){\footnotesize $\Pi$}
\put(160,40){\footnotesize $\Gamma$}
\put(47,40){\footnotesize $\Gamma$}
\put(90,45){\footnotesize $\Delta u=1$}
\put(95,75){\footnotesize $\Omega(u)$}
\put(46,90){\footnotesize $B_1^+$}
\put(22,27){\footnotesize $u=0$}
\put(22,17){\footnotesize $|\nabla u|=0$}
\put(165,27){\footnotesize $u=0$}
\put(160,17){\footnotesize $|\nabla u|=0$}
\put(100,0){\footnotesize $u=0$}
\end{picture}
\caption{Touch between the free boundary $\Gamma=\partial\Omega(u)$
and the fixed boundary $\Pi$ in problem \eqref{eq:tangential}.}
\label{fig:tangential}
\end{figure}
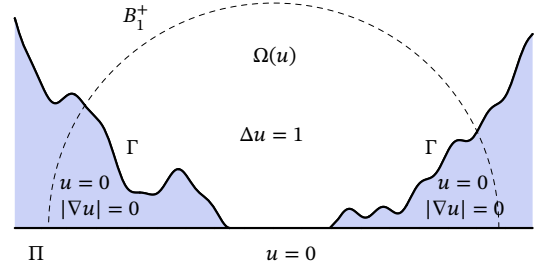

\subsubsection{Touch between free and fixed boundary}
\label{sec:touch}

In \cite{MR1359745} (joint with one of us) and her follow-up paper  \cite{MR1392033}, 
Uraltseva studied the obstacle problem close to a Dirichlet data, for
smooth boundaries, where she proves that the free  boundary touches
the fixed boundary tangentially. The idea  seemed to be  inspired by
related works with Oliker (see Section~\ref{sec:geom-eq}) and  the Dam-problem in filtration.

During the potential theory program at  Institute Mittag-Leffler (1999--2000) she started working on free boundary problems that originated in potential theory. Specifically, the harmonic continuation problem in potential theory, that was strongly tied to obstacle problem, but with the lack of having a sign for the solution function. The simplest way to  formulate this problem is as follows:
\begin{equation}\label{eq:tangential}
\begin{aligned}
&\Delta u =\chi_{\Omega(u)}\quad\text{in } B_1^+,\\ &\qquad\qquad\text{with }\Omega(u)\coloneq\{u=|\nabla u|=0\}^c\\
&u=0 \quad \text{on }\Pi\cap B_1,
\end{aligned}
\end{equation}
where $B_1^+= \{|x| < 1, \ x_1 > 0\}$  and $\Pi=\{x_1=0\}$;  see Figure~\ref{fig:tangential}.
The  question of interest was the behavior of the free boundary $\Gamma=\partial\Omega(u)$ close to the fixed boundary $\Pi$.

In  \cite{MR1810607},  and several follow-up papers in the parabolic regime, she shows that the free boundary $\Gamma$   is a graph of a $C^1$-function close to points on $\Pi$, where $\Gamma \cap B_1^+$ touches $\Pi$, or comes too close to $\Pi$.
 
To   prove   this, and the related parabolic results, there was a need for developing new tools and approaches. This was possible partly due to the  availability of monotonicity formulas, such as that of Alt, Caffarelli, and Friedman \cite{MR732100}. One version of the latter  asserts that  for continuous subharmonic  functions $h_1,$ $h_2$   in $B_R(x^0)$,
 satisfying  $h_1h_2 =0$, and $h_1(x^0) = h_2 (x^0)  =0$, we have 
 $\varphi (r) \nearrow $ for $0<r< R$,   where 
\begin{equation}\label{monot-formula}
\varphi (r)=\phi (r,h_1,x^0)\,\phi (r,h_2,x^0)
\end{equation}
with 
$$
\phi (r,h_i,x^0)\coloneq \frac{ 1}{ r^{2}} \int_{B(x^0,r)}
 \frac{| \nabla   h_i |^2 dx}{|x-x^0|^{n-2}} .
$$
One can use the monotonicity of the function $\varphi(r)$ to prove several important properties for $u$ and the free boundary. Indeed, one first  extends $u$ to be zero in $B_1^-= \{|x| < 1, \ x_1 < 0\}$ and  applies the monotonicity formula \eqref{monot-formula} to $h_{1}=(\partial_e u)^+$ and $h_{2}=(\partial_e u)^-$, where $e$ is any vector tangent to the plane $\{x_1=0\}$. Using the fact that at least one of the sets $\{ \pm \partial_e u > 0\}$  has positive volume density at $x^0$, we shall have 
$$c_0| \nabla \partial_e u (x^0)|^4 = \lim_{r\to 0}  \varphi (r) \leq  \varphi (1) \leq C_0.
$$
Combining this with  equation \eqref{eq:tangential} we obtain the bound for $u_{x_1x_1} (x^0)$.
From here, the  uniform $C^{1,1}$ regularity 
for $u$ in $B_{1/2}^+$ follows. 

The $C^{1,1}$ regularity is instrumental for any analysis of the properties of the free boundary. Indeed, to study the free boundary at points where it touches the fixed boundary, one needs to rescale the solution quadratically,
$u_r (x) = u(rx + x^0) /r^2$, which keeps the equation invariant. Indeed, this scaling and ``blow-up''\footnote{Blow-up refers to  $\lim_{r \to 0} u_r (x)$, whenever it exists. }
brings one to a global setting of equation  \eqref{eq:tangential}  in $\R^n_+$, where solutions can be  classified (in a rotated system) as one of the following:
\begin{itemize}
\item[(i)] $u(x)=\frac12 x_1^2 + ax_1x_2 + \alpha x_1$\quad $(a>0,\ \alpha\in\R)$,
\item[(ii)] $u(x)= \frac12 ((x_1 - a)^+)^2$ \quad $(a>0)$.
\end{itemize}
The proof of the classification of global solutions uses an array of geometric tools and the monotonicity function $\varphi (r)$,  implying that if $\{u =0 \}\cap \{x_1 > 0\} \neq \emptyset$, then $\partial_e u \equiv 0$, for any direction $e$ tangential to $\Pi$. The case when this set is empty is easily handled by Liouville's theorem.

Once  this classification is done, one can argue by  indirect methods that the free boundary  $\partial \{u > 0\} \cap \{x_1 >0\}$ approaches the fixed one, at touching points, tangentially,  and that it is a $C^1$-graph locally, which is optimal in the sense that  (in general) it cannot be $C^{1,{\rm Dini}}$.

\subsubsection{Two-phase obstacle type problems}

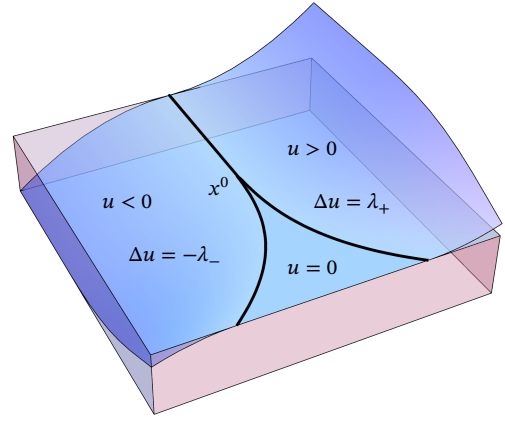
\begin{figure}
  \centering
\begin{picture}(200,160)
\put(0,0){\includegraphics[height=160pt]{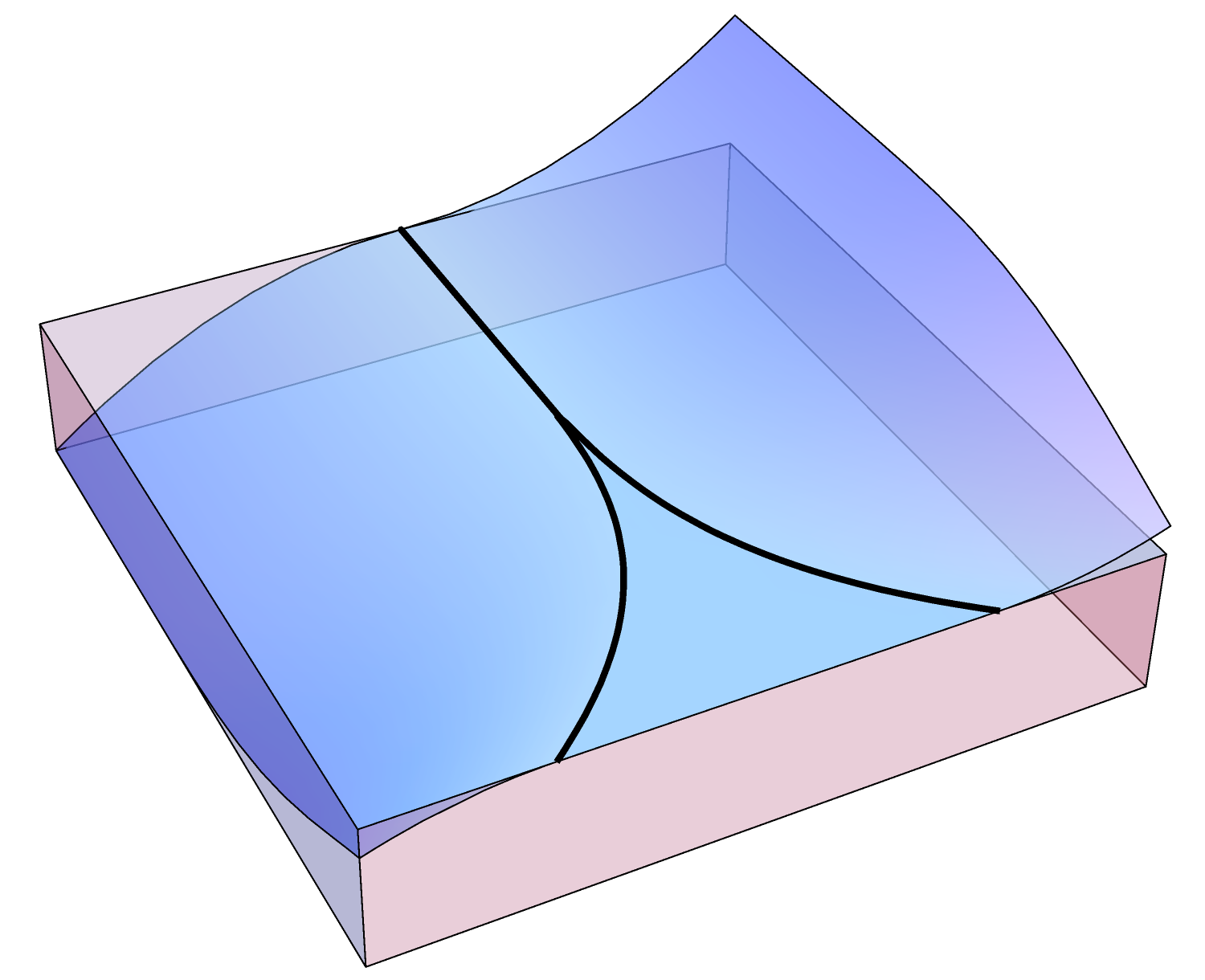}}
\put(40,80){\footnotesize $u<0$}
\put(50,60){\footnotesize $\Delta u=-\lambda_-$}
\put(110,100){\footnotesize $u>0$}
\put(120,80){\footnotesize $\Delta u=\lambda_+$}
\put(110,55){\footnotesize $u=0$}
\put(80,85){\footnotesize $x^0$}
\end{picture}
\caption{Two-phase problem: branch point $x^0$.}
\label{fig:two-phase}
\end{figure}

If one  considers extension of  equation \eqref{eq:tangential}  into $B_1$,  by an odd reflection, then one obtains a specific example of a general problem that is referred to as two-phase obstacle problem, and is formulated as 
$$\Delta u = \lambda_+ \chi_{\{u>0\}} -
\lambda_- \chi_{\{u<0\}}\quad\text{in } B_1(0),$$
where $\lambda_\pm$ are positive bounded Lipschitz functions.  Figure~\ref{fig:two-phase} illustrates this problem.

In \cite{MR2340105},  Nina Uraltseva (with co-authors) proves that at any  branch point $x^0 \in \partial \{u>0\} \cap \partial \{u<0\}$ with $u(x^0)=  |\nabla u(x^0)| =0 $, the free boundaries  $\partial \{ u>0\}\cap B_{r_0}(0)$ and $\partial \{ u<0\}\cap B_{r_0}(0)$ are $C^{1}$-surfaces, that touch each other tangentially at $x^0$.

The proof of this and several similar results (also in parabolic setting) relies heavily on the  monotonicity function $\varphi$ mentioned above as well as  on the balanced energy functional 
\begin{multline}\label{weiss}
 \Phi_{x_0}(r) \coloneq r^{-n-2} \int_{B_r(x_0)} \left( 
{\vert \nabla u \vert}^2 +  \lambda_+ u^+
+ \lambda_- u^-\right)\\   
- 2 r^{-n-3}  \int_{\partial B_r(x_0)}
u^2  ,
\end{multline}
which is strictly monotone in $r$, unless $u$ is homogeneous.
Using these two monotonicity functionals in combination with geometric tools bring us to the fact that any global solution $u_0$ to the two-phase problem is one-dimensional and, in a rotated and translated system of coordinates,  
$$u_0= \frac{\lambda_+}{2} (x_1^+)^2 - \frac{\lambda_-}{2} (x_1^-)^2.$$
From here one uses a revised form of the so-called directional monotonicity argument  of Luis Caffarelli, that in this setting boils down to the fact that close to branch points $x^0$ one can show 
that in a suitable cone of directions $ \mathcal{C}$ one has
$\partial_e u \geq 0$,  in $B_r (x^0)$ for $e \in  \mathcal{C}$  and $r$ universal. This in particular implies that the free boundaries
 $\partial \{ \pm u >0 \} $ 
 are Lipschitz graphs locally close to branch points.

The approaches here generated further application of the techniques to  problems with hysteresis \cite{MR3352792,MR3393312}.

\subsubsection{Free boundaries for weakly coupled  systems}

In her work with coauthors \cite{MR3350233}, Uraltseva considers the following vectorial   energy  minimizing functional:
$$
E(\mathbf{u})\coloneq \int_{B_1} \left(|\nabla \mathbf{u}|^2+ 2 |\mathbf{u}|  \right) dx .
$$
Here $B_1$ is the unit ball in $\R^n$ ($n\geq 1$),  and we minimize over the Sobolev space
 $\mathbf{g} + W_{0}^{1,2} (B_1; \R^N)$  for some smooth boundary values  $\mathbf{g}=(g_1, \dots, g_N)$.
The minimizer(s)  are   vector-valued functions $\mathbf{u} =(u_1, \dots , u_N)$,  with components $u_i$ satisfying
$$
\Delta u_i=\frac{u_i}{|\mathbf{u}|}, \quad i=1,\dots, N.
$$
Since the set  $\{ |\mathbf{u}| > 0\}$  competes with the Dirichlet energy, by taking  the boundary values  small we may obtain   $\{ \mathbf{u} = \mathbf{0} \} \neq \emptyset$, 
which is in contrast   to standard variational problems.
The set $\partial \{  |\mathbf{u}| > 0 \} $ is called the free boundary. 
One observes that when $N=1$ (scalar case) then we fall back to the two-phase problem.

Simple examples of solutions to this problem are
\begin{itemize}
\item[(i)] $u_i = \alpha_i P(x)$, with $P(x) \geq 0$, $\Delta P (x)=1$, and $\sum_{i=1}^N \alpha_i^2 = 1$,\vspace{1ex} 
\item[(ii)] $u_i = \frac{\alpha_i}2 (x_1^+)^2 + \frac{\beta_i}2 (x_1^-)^2$  \quad (2-phase),\\
$\sum_{i=1}^N \alpha_i^2 = 1$,\quad  $\sum_{i=1}^N \beta_i^2 = 1$,\vspace{1ex}   
\item[(iii)] $u_i = \frac{\alpha_i}2(x_1^+)^2$ \quad (1-phase),\\
$\sum_{i=1}^N \alpha_i^2 = 1$.
\end{itemize}
Using the vectorial version of the  monotonicity formula \eqref{weiss}, one can show that $\mathbf{u}$ has a quadratic growth away from the free boundary.

The regularity of the free boundary follows through the homogeneity improvement approach with the so-called epiperimetric inequality, which is used to show that the functional 
$$
\mathcal{M}(\mathbf{v})\coloneq \int_{B_1(0)} (|\nabla \mathbf{v}|^2 + 2 |\mathbf{v}|) -  \int_{\partial B_1(0)} |\mathbf{v}|^2
$$
satisfies
$$
| \mathcal{M}(\mathbf{u}_{r_1}) - \mathcal{M}(\mathbf{u}_{r_2})| \leq 
c |r_2 -r_1|^\alpha,\quad \alpha>0,
$$
where $\mathbf{u}_r = \mathbf{u} (x^0+rx)/r^2$, and $x^0 $ is such that $\mathbf{u}_r$ is close to rotated version of a half-space solution of type 
$\mathbf{h}=\frac12 (x_1^+)^2 \mathbf{e}$.

This, in particular,  gives uniqueness of the blow-ups, and can be used to show that (in a rotated system of coordinates) there exist  $\beta' >0$,  $r_0>0$, and $C<\infty$ such that
\begin{multline*}
\int_{\partial B_1(0)} \left\vert \mathbf{u}_r - \mathbf{h}
\right\vert \leq Cr^{\beta'}\\
\text{for every } x^0\in 
\mathcal{R}_\mathbf{u }\text{ and every } r \le r_0,
\end{multline*}
where $\mathcal{R}_\mathbf{u}$ is the set of  free boundary points whose blow-ups are half-spaces.
This implies that $\mathcal{R}_\mathbf{u}$   is  locally in $B_{1/2}$ a $C^{1,\beta}$-surface.

\begin{figure}
\centering
\includegraphics[height=\columnwidth]{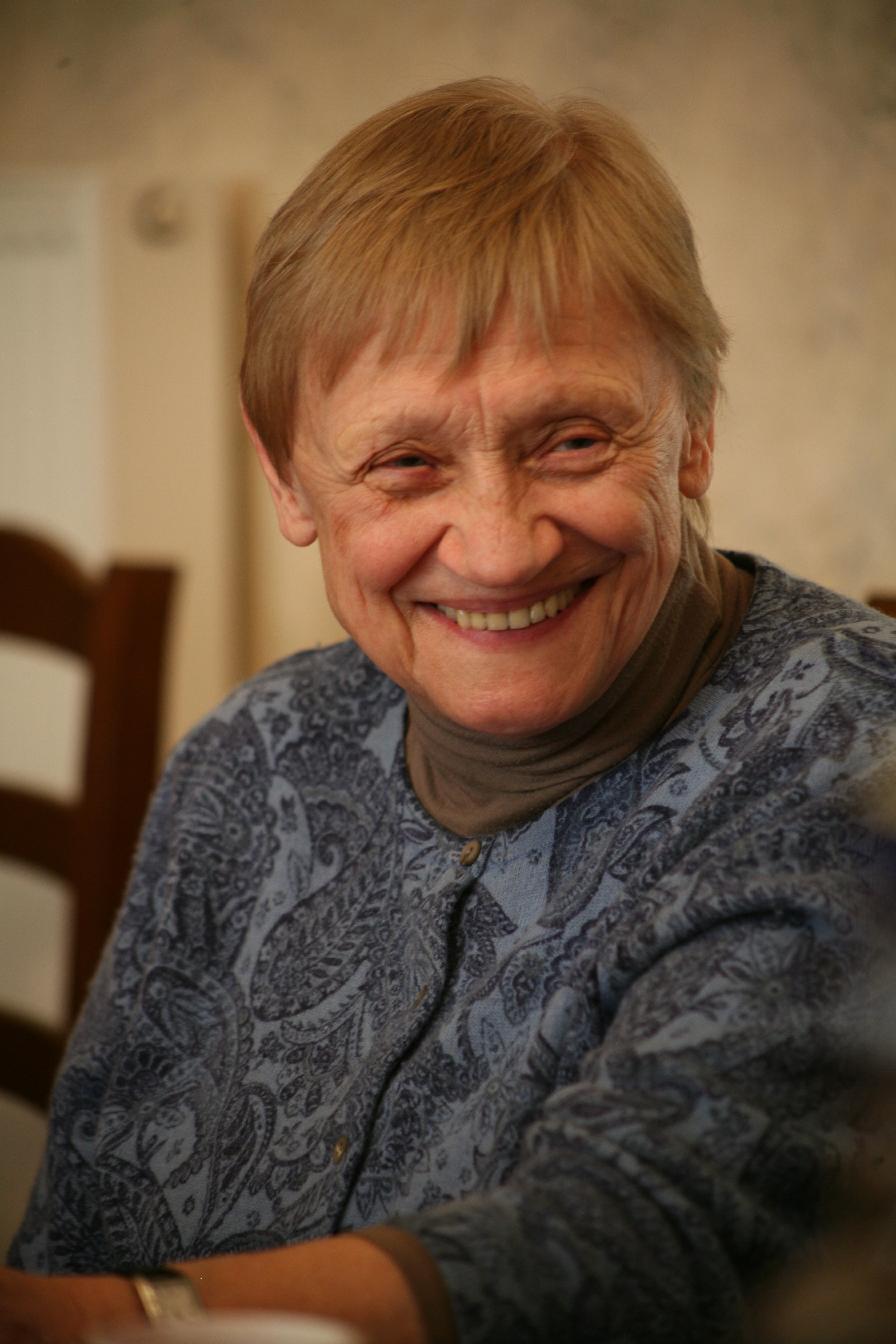}
\caption{Nina Uraltseva in 2013.}
\label{fig:2013}
\end{figure}

\section*{Nina's Impact}

Nina Uraltseva  has over 100 publications\footnote{See:  \url{https://www.scilag.net/profile/nina-uraltseva}.} and over 8000 citations in MathSciNet. Her famous book \emph{Linear and Quasilinear Equations of Parabolic Type} \cite{MR0241821,MR0241822} (joint with Ladyzhenskaya and Solonnikov) has over 4600 citation, and the elliptic version of this book  \cite{MR0211073,MR0244627} (joint with Ladyzhenskaya), has over 1600 citation in MathSciNet. This naturally gives a picture of a mathematician with tremendous impact on the filed of partial differential equations. 
Needless to say that, even though there are many new  books on the topic of PDEs, 
these books stay equally important and extremely valuable  to many PhD students and  early career analysts.

Nina Uraltseva has, over the years, contributed  to the mathematical community by serving on many important committees; e.g.,\ chairing  the PDE Panel of the International Congress of Mathematicians in Berlin, Germany, 1998, and  the Prize Committee of the European Congress of Mathematics in Stockholm, Sweden, 2004. She also served as an expert for research foundations such as the European Research Council and the Russian Foundation for Basic Research.

She has been acting as an editor for several journals,\footnote{Editor in Chief for \emph{Proceedings of St.\ Petersburg Math.\ Society} and \emph{Journal of Problems in Mathematical Analysis};  member of the editorial committee for \emph{Algebra and Analysis} (translated in \emph{St. Petersburg Mathematical Journal}), \emph{Vestnik St.\ Petersburg State University}, \emph{Lithuanian Mathematical Journal}.} and has been a frequent visitor of many universities all over the world and presented talks at various international conferences and schools. In her role as a world leading  expert in analysis of PDEs she  has captured the attention of many female  students in all areas of mathematics, and attracted them to further pursue research and start a career in mathematics. Her motivational talks at many conferences, specially meetings related to ``connection to women,'' have been an important factor in attracting several females to mathematics.

The instructional aspect of her work and her dedication to educating PhD students,\footnote{Uraltseva has supervised 13 PhD students, four of which have habilitated.} as well as unselfishly being available to students and colleagues for discussions and brain-storming of their problems, make her one of the most prominent and devoted persons to the mathematical community.

Nina Uraltseva has dedicated her life to mathematics, and in her scientific journey through the years she has made many friends all over the world. Her kind personality and utmost politeness on one side and her unbiased style and open mindedness towards diverse mathematical problems have made  her extremely popular  among  colleagues and  students, and  not only as a mathematician but also as a human being.

\bibliographystyle{amsplain}

\let\OLDthebibliography\thebibliography
\renewcommand\thebibliography[1]{
  \OLDthebibliography{#1}
  \setlength{\parskip}{0pt}
  \setlength{\itemsep}{0pt plus 0.3ex}
}

{\small \bibliography{Uraltseva-AMS}}

\providecommand{\bysame}{\leavevmode\hbox to3em{\hrulefill}\thinspace}
\providecommand{\MR}{\relax\ifhmode\unskip\space\fi MR }
\providecommand{\MRhref}[2]{%
  \href{http://www.ams.org/mathscinet-getitem?mr=#1}{#2}
}
\providecommand{\href}[2]{#2}
\begin{thebibliography}{10}

\bibitem{MR732100}
Hans~Wilhelm Alt, Luis~A. Caffarelli, and Avner Friedman, \emph{Variational
  problems with two phases and their free boundaries}, Trans. Amer. Math. Soc.
  \textbf{282} (1984), no.~2, 431--461. \MR{732100}

\bibitem{MR3350233}
John Andersson, Henrik Shahgholian, Nina~N. Uraltseva, and Georg~S. Weiss,
  \emph{Equilibrium points of a singular cooperative system with free
  boundary}, Adv. Math. \textbf{280} (2015), 743--771. \MR{3350233}

\bibitem{MR1810607}
D.~E. Apushkinskaya, H.~Shahgholian, and N.~N. Uraltseva, \emph{Boundary
  estimates for solutions of a parabolic free boundary problem}, Zap. Nauchn.
  Sem. S.-Peterburg. Otdel. Mat. Inst. Steklov. (POMI) \textbf{271} (2000),
  no.~Kraev. Zadachi Mat. Fiz. i Smezh. Vopr. Teor. Funkts. 31, 39--55, 313.
  \MR{1810607}

\bibitem{MR1359745}
D.~E. Apushkinskaya and N.~N. Ural\cprime~tseva, \emph{On the behavior of the
  free boundary near the boundary of the domain}, Zap. Nauchn. Sem.
  S.-Peterburg. Otdel. Mat. Inst. Steklov. (POMI) \textbf{221} (1995),
  no.~Kraev. Zadachi Mat. Fiz. i Smezh. Voprosy Teor. Funktsi\u{\i}. 26, 5--19,
  253. \MR{1359745}

\bibitem{MR3393312}
D.~E. Apushkinskaya and N.~N. Uraltseva, \emph{Free boundaries in problems with
  hysteresis}, Philos. Trans. Roy. Soc. A \textbf{373} (2015), no.~2050,
  20140271, 10. \MR{3393312}

\bibitem{MR3352792}
Darya~E. Apushkinskaya and Nina~N. Uraltseva, \emph{On regularity properties of
  solutions to the hysteresis-type problem}, Interfaces Free Bound. \textbf{17}
  (2015), no.~1, 93--115. \MR{3352792}

\bibitem{MR1374941}
A.~Arkhipova and N.~Uraltseva, \emph{Sharp estimates for solutions of a
  parabolic {S}ignorini problem}, Math. Nachr. \textbf{177} (1996), 11--29.
  \MR{1374941}

\bibitem{MR869237}
A.~A. Arkhipova and N.~N. Ural\cprime~tseva, \emph{Regularity of solutions of
  diagonal elliptic systems under convex constraints on the boundary of the
  domain}, Zap. Nauchn. Sem. Leningrad. Otdel. Mat. Inst. Steklov. (LOMI)
  \textbf{152} (1986), no.~Kraev. Zadachi Mat. Fiz. i Smezhnye Vopr. Teor.
  Funktsi\u{\i}18, 5--17, 181. \MR{869237}

\bibitem{MR841488}
\bysame, \emph{Regularity of the solution of a problem with a two-sided
  constraint on the boundary}, Vestnik Leningrad. Univ. Mat. Mekh. Astronom.
  (1986), no.~vyp. 1, 3--10, 133. \MR{841488}

\bibitem{MR928154}
\bysame, \emph{Regularity of the solutions of variational inequalities with
  convex constraints on the boundary of the domain for nonlinear operators with
  a diagonal principal part}, Vestnik Leningrad. Univ. Mat. Mekh. Astronom.
  (1987), no.~vyp. 3, 13--19, 127. \MR{928154}

\bibitem{MR964910}
\bysame, \emph{Regularity of the solution of a problem with a two-sided limit
  on a boundary for elliptic and parabolic equations}, Trudy Mat. Inst.
  Steklov. \textbf{179} (1988), 5--22, 241, Translated in Proc. Steklov Inst.
  Math. {\bf 1989}, no. 2, 1--19, Boundary value problems of mathematical
  physics, 13 (Russian). \MR{964910}

\bibitem{MR542512}
L.~A. Caffarelli, \emph{Further regularity for the {S}ignorini problem}, Comm.
  Partial Differential Equations \textbf{4} (1979), no.~9, 1067--1075.
  \MR{542512}

\bibitem{MR709038}
E.~DiBenedetto, \emph{{$C\sp{1+\alpha }$} local regularity of weak solutions of
  degenerate elliptic equations}, Nonlinear Anal. \textbf{7} (1983), no.~8,
  827--850. \MR{709038}

\bibitem{MR672713}
Lawrence~C. Evans, \emph{A new proof of local {$C\sp{1,\alpha }$} regularity
  for solutions of certain degenerate elliptic p.d.e}, J. Differential
  Equations \textbf{45} (1982), no.~3, 356--373. \MR{672713}

\bibitem{MR983300}
Gerhard Huisken, \emph{Nonparametric mean curvature evolution with boundary
  conditions}, J. Differential Equations \textbf{77} (1989), no.~2, 369--378.
  \MR{983300}

\bibitem{MR620584}
David Kinderlehrer, \emph{The smoothness of the solution of the boundary
  obstacle problem}, J. Math. Pures Appl. (9) \textbf{60} (1981), no.~2,
  193--212. \MR{620584}

\bibitem{MR0241821}
O.~A. Lady\v{z}enskaja, V.~A. Solonnikov, and N.~N. Ural\cprime~ceva,
  \emph{{\cyr Line\u{i}nye i kvaziline\u{i}nye uravneniya parabolicheskogo
  tipa}}, Izdat. ``Nauka'', Moscow, 1967. \MR{0241821}

\bibitem{MR0241822}
\bysame, \emph{Linear and quasilinear equations of parabolic type},
  Translations of Mathematical Monographs, Vol. 23, American Mathematical
  Society, Providence, R.I., 1968, Translated from the Russian by S. Smith.
  \MR{0241822}

\bibitem{MR0141891}
O.~A. Lady\v{z}enskaja and N.~N. Ural\cprime~ceva, \emph{A boundary-value
  problem for linear and quasi-linear parabolic equations}, Dokl. Akad. Nauk
  SSSR \textbf{139} (1961), 544--547. \MR{0141891}

\bibitem{MR0141874}
\bysame, \emph{Differential properties of bounded generalized solutions of
  multidimensional quasilinear elliptic equations and variational problems},
  Dokl. Akad. Nauk SSSR \textbf{138} (1961), 29--32. \MR{0141874}

\bibitem{MR0149075}
\bysame, \emph{Quasilinear elliptic equations and variational problems in
  several independent variables}, Uspehi Mat. Nauk \textbf{16} (1961), no.~1
  (97), 19--90. \MR{0149075}

\bibitem{MR0150447}
\bysame, \emph{Regularity of generalized solutions of quasilinear elliptic
  equations}, Dokl. Akad. Nauk SSSR \textbf{140} (1961), 45--47. \MR{0150447}

\bibitem{MR0181837}
\bysame, \emph{A boundary-value problem for linear and quasi-linear parabolic
  equations. {I}, {II}, {III}}, Iaz. Akad. Nauk SSSR Ser. Mat. 26 (1962), 5-52;
  ibid. 26 (1962), 753- 780; ibid. \textbf{27} (1962), 161--240. \MR{0181837}

\bibitem{MR0147786}
\bysame, \emph{The first boundary-value problem for quasi-linear second-order
  parabolic equations of general type}, Dokl. Akad. Nauk SSSR \textbf{147}
  (1962), 28--30. \MR{0147786}

\bibitem{MR0211073}
\bysame, \emph{{\cyr Line\u{i}nye i kvaziline\u{i}nye uravneniya
  \`ellipticheskogo tipa}}, Izdat. ``Nauka'', Moscow, 1964. \MR{0211073}

\bibitem{MR149076}
O.~A. Ladyzhenskaia and N.~N. Ural\cprime~tzeva, \emph{On the smoothness of
  weak solutions of quasilinear equations in several variables and of
  variational problems}, Comm. Pure Appl. Math. \textbf{14} (1961), 481--495.
  \MR{149076}

\bibitem{MR265745}
O.~A. Ladyzhenskaya and N.~N. Ural\cprime~tseva, \emph{Local estimates for
  gradients of solutions of non-uniformly elliptic and parabolic equations},
  Comm. Pure Appl. Math. \textbf{23} (1970), 677--703. \MR{265745}

\bibitem{MR0509265}
\bysame, \emph{{\cyr Line\u{i}nye i kvaziline\u{i}nye uravneniya} \`e{\cyr
  llipticheskogo tipa}}, Izdat. ``Nauka'', Moscow, 1973, Second edition,
  revised. \MR{0509265}

\bibitem{MR878325}
\bysame, \emph{A survey of results on the solvability of boundary value
  problems for uniformly elliptic and parabolic second-order quasilinear
  equations having unbounded singularities}, Uspekhi Mat. Nauk \textbf{41}
  (1986), no.~5(251), 59--83, 262. \MR{878325}

\bibitem{MR0244627}
Olga~A. Ladyzhenskaya and Nina~N. Ural'tseva, \emph{Linear and quasilinear
  elliptic equations}, Academic Press, New York-London, 1968, Translated from
  the Russian by Scripta Technica, Inc, Translation editor: Leon Ehrenpreis.
  \MR{0244627}

\bibitem{MR721568}
John~L. Lewis, \emph{Regularity of the derivatives of solutions to certain
  degenerate elliptic equations}, Indiana Univ. Math. J. \textbf{32} (1983),
  no.~6, 849--858. \MR{721568}

\bibitem{MR868523}
Paolo Marcellini, \emph{On the definition and the lower semicontinuity of
  certain quasiconvex integrals}, Ann. Inst. H. Poincar\'{e} Anal. Non
  Lin\'{e}aire \textbf{3} (1986), no.~5, 391--409. \MR{868523}

\bibitem{MR969900}
\bysame, \emph{Regularity of minimizers of integrals of the calculus of
  variations with nonstandard growth conditions}, Arch. Rational Mech. Anal.
  \textbf{105} (1989), no.~3, 267--284. \MR{969900}

\bibitem{MR1094446}
\bysame, \emph{Regularity and existence of solutions of elliptic equations with
  {$p,q$}-growth conditions}, J. Differential Equations \textbf{90} (1991),
  no.~1, 1--30. \MR{1094446}

\bibitem{MR4258810}
Giuseppe Mingione and Vicen\c{t}iu R\v{a}dulescu, \emph{Recent developments in
  problems with nonstandard growth and nonuniform ellipticity}, J. Math. Anal.
  Appl. \textbf{501} (2021), no.~1, Paper No. 125197, 41. \MR{4258810}

\bibitem{MR821477}
A.~I. Nazarov and N.~N. Ural\cprime~tseva, \emph{Convex-monotone hulls and an
  estimate of the maximum of the solution of a parabolic equation}, Zap.
  Nauchn. Sem. Leningrad. Otdel. Mat. Inst. Steklov. (LOMI) \textbf{147}
  (1985), 95--109, 204--205, Boundary value problems of mathematical physics
  and related problems in the theory of functions, No. 17. \MR{821477}

\bibitem{MR2760150}
\bysame, \emph{The {H}arnack inequality and related properties of solutions of
  elliptic and parabolic equations with divergence-free lower-order
  coefficients}, Algebra i Analiz \textbf{23} (2011), no.~1, 136--168.
  \MR{2760150}

\bibitem{MR1334142}
V.~I. Oliker and N.~N. Ural\cprime~tseva, \emph{Long time behavior of flows
  moving by mean curvature}, Nonlinear evolution equations, Amer. Math. Soc.
  Transl. Ser. 2, vol. 164, Amer. Math. Soc., Providence, RI, 1995,
  pp.~163--170. \MR{1334142}

\bibitem{MR1483640}
Vladimir~I. Oliker and Nina~N. Ural\cprime~tseva, \emph{Long time behavior of
  flows moving by mean curvature. {II}}, Topol. Methods Nonlinear Anal.
  \textbf{9} (1997), no.~1, 17--28. \MR{1483640}

\bibitem{MR1193345}
Vladimir~I. Oliker and Nina~N. Uraltseva, \emph{Evolution of nonparametric
  surfaces with speed depending on curvature. {II}. {T}he mean curvature case},
  Comm. Pure Appl. Math. \textbf{46} (1993), no.~1, 97--135. \MR{1193345}

\bibitem{MR1246345}
\bysame, \emph{Evolution of nonparametric surfaces with speed depending on
  curvature. {III}. {S}ome remarks on mean curvature and anisotropic flows},
  Degenerate diffusions ({M}inneapolis, {MN}, 1991), IMA Vol. Math. Appl.,
  vol.~47, Springer, New York, 1993, pp.~141--156. \MR{1246345}

\bibitem{MR2962060}
Arshak Petrosyan, Henrik Shahgholian, and Nina Uraltseva, \emph{Regularity of
  free boundaries in obstacle-type problems}, Graduate Studies in Mathematics,
  vol. 136, American Mathematical Society, Providence, RI, 2012. \MR{2962060}

\bibitem{MR986184}
Rainer Schumann, \emph{Regularity for {S}ignorini's problem in linear
  elasticity}, Manuscripta Math. \textbf{63} (1989), no.~3, 255--291.
  \MR{986184}

\bibitem{MR1950478}
Henrik Shahgholian and Nina Uraltseva, \emph{Regularity properties of a free
  boundary near contact points with the fixed boundary}, Duke Math. J.
  \textbf{116} (2003), no.~1, 1--34. \MR{1950478}

\bibitem{MR2340105}
Henrik Shahgholian, Nina Uraltseva, and Georg~S. Weiss, \emph{The two-phase
  membrane problem---regularity of the free boundaries in higher dimensions},
  Int. Math. Res. Not. IMRN (2007), no.~8, Art. ID rnm026, 16. \MR{2340105}

\bibitem{MR727034}
Peter Tolksdorf, \emph{Regularity for a more general class of quasilinear
  elliptic equations}, J. Differential Equations \textbf{51} (1984), no.~1,
  126--150. \MR{727034}

\bibitem{MR369884}
Neil~S. Trudinger, \emph{Linear elliptic operators with measurable
  coefficients}, Ann. Scuola Norm. Sup. Pisa Cl. Sci. (3) \textbf{27} (1973),
  265--308. \MR{369884}

\bibitem{MR474389}
K.~Uhlenbeck, \emph{Regularity for a class of non-linear elliptic systems},
  Acta Math. \textbf{138} (1977), no.~3-4, 219--240. \MR{474389}

\bibitem{MR0126742}
N.~N. Ural\cprime~ceva, \emph{Regularity of solutions of multidimensional
  elliptic equations and variational problems}, Soviet Math. Dokl. \textbf{1}
  (1960), 161--164. \MR{0126742}

\bibitem{MR0142886}
\bysame, \emph{Boundary-value problems for quasi-linear elliptic equations and
  systems with principal part of divergence type}, Dokl. Akad. Nauk SSSR
  \textbf{147} (1962), 313--316. \MR{0142886}

\bibitem{MR0140817}
\bysame, \emph{General second-order quasi-linear equations and certain classes
  of systems of equations of elliptic type}, Dokl. Akad. Nauk SSSR \textbf{146}
  (1962), 778--781. \MR{0140817}

\bibitem{MR0244628}
\bysame, \emph{Degenerate quasilinear elliptic systems}, Zap. Nau\v{c}n. Sem.
  Leningrad. Otdel. Mat. Inst. Steklov. (LOMI) \textbf{7} (1968), 184--222.
  \MR{0244628}

\bibitem{MR0364860}
\bysame, \emph{Nonlinear boundary value problems for equations of minimal
  surface type}, Trudy Mat. Inst. Steklov. \textbf{116} (1971), 217--226, 237,
  Boundary value problems of mathematical physics, 7. \MR{0364860}

\bibitem{MR0313623}
\bysame, \emph{The regularity of the solutions of variational inequalities},
  Zap. Nau\v{c}n. Sem. Leningrad. Otdel. Mat. Inst. Steklov. (LOMI) \textbf{27}
  (1972), 211--219, Boundary value problems of mathematical physics and related
  questions in the theory of functions, 6. \MR{0313623}

\bibitem{MR0638359}
\bysame, \emph{The solvability of the capillarity problem}, Vestnik Leningrad.
  Univ. (1973), no.~19 Mat. Meh. Astronom. Vyp. 4, 54--64, 152. \MR{0638359}

\bibitem{MR0638360}
\bysame, \emph{The solvability of the capillarity problem. {II}}, Vestnik
  Leningrad. Univ. (1975), no.~1, Mat. Meh. Astronom. vyp. 1, 143--149, 191,
  Collection of articles dedicated to the memory of Academician V. I. Smirnov.
  \MR{0638360}

\bibitem{MR508511}
\bysame, \emph{Strong solutions of the generalized {S}ignorini problem},
  Sibirsk. Mat. Zh. \textbf{19} (1978), no.~5, 1204--1212, 1216. \MR{508511}

\bibitem{MR660089}
N.~N. Ural\cprime~tseva, \emph{Estimates of the maximum moduli of gradients for
  solutions of capillarity problems}, Zap. Nauchn. Sem. Leningrad. Otdel. Mat.
  Inst. Steklov. (LOMI) \textbf{115} (1982), 274--284, 312, Boundary value
  problems of mathematical physics and related questions in the theory of
  functions, 14. \MR{660089}

\bibitem{MR775926}
\bysame, \emph{H\"{o}lder continuity of gradients of solutions of parabolic
  equations with boundary conditions of {S}ignorini type}, Dokl. Akad. Nauk
  SSSR \textbf{280} (1985), no.~3, 563--565. \MR{775926}

\bibitem{MR860572}
\bysame, \emph{Estimation on the boundary of the domain of derivatives of
  solutions of variational inequalities}, Linear and nonlinear boundary value
  problems. {S}pectral theory ({R}ussian), Probl. Mat. Anal., vol.~10,
  Leningrad. Univ., Leningrad, 1986, Translated in J. Soviet Math. {\bf 45}
  (1989), no. 3, 1181--1191, pp.~92--105, 213. \MR{860572}

\bibitem{MR934318}
\bysame, \emph{Estimates of derivatives of solutions of elliptic and parabolic
  inequalities}, Proceedings of the {I}nternational {C}ongress of
  {M}athematicians, {V}ol. 1, 2 ({B}erkeley, {C}alif., 1986), Amer. Math. Soc.,
  Providence, RI, 1987, pp.~1143--1149. \MR{934318}

\bibitem{MR933999}
\bysame, \emph{On the regularity of solutions of variational inequalities},
  Uspekhi Mat. Nauk \textbf{42} (1987), no.~6(258), 151--174, 248. \MR{933999}

\bibitem{MR1392033}
\bysame, \emph{{$C^1$} regularity of the boundary of a noncoincident set in a
  problem with an obstacle}, Algebra i Analiz \textbf{8} (1996), no.~2,
  205--221. \MR{1392033}

\bibitem{MR725829}
N.~N. Ural\cprime~tseva and A.~B. Urdaletova, \emph{Boundedness of gradients of
  generalized solutions of degenerate nonuniformly elliptic quasilinear
  equations}, Vestnik Leningrad. Univ. Mat. Mekh. Astronom. (1983), no.~vyp. 4,
  50--56. \MR{725829}

\end{thebibliography}

\begin{small}
\subsubsection*{Credits}

\begin{description}\itemsep=0pt
\item Opening photo is courtesy of Norayr Matevosyan.
\item Figures~\ref{fig:1951} and \ref{fig:seminar} are from Nina Uraltseva's personal archive (Figure~\ref{fig:seminar} photograph by Nina Alovert).
\item Figures~\ref{fig:bdry-obst}--\ref{fig:two-phase} are courtesy of Arshak Petrosyan.
\item Figure~\ref{fig:2013} is courtesy of Sophia Nazarova.
\end{description}
\end{small}

\end{document}